\documentclass{amsart}

\def\Z{\mathbf{Z}}
\def\C{\mathbf{C}}

\DeclareFontFamily{OT1}{rsfs}{}
\DeclareFontShape{OT1}{rsfs}{n}{it}{<-> rsfs10}{}
\DeclareMathAlphabet{\mathscr}{OT1}{rsfs}{n}{it}

\newcommand{\et}{\'{e}t}
\newcommand{\rightrightarrows}{\rightarrow}

\newcommand{\bigdot}{\bullet}

\DeclareMathOperator{\sHom}{\mathscr{H}om}
\DeclareMathOperator{\Coh}{Coh}

\DeclareMathOperator{\an}{an}

\DeclareMathOperator{\Mod}{\mathcal{M}}

\DeclareMathOperator{\GL}{GL}

\DeclareMathOperator{\Ext}{Ext}

\DeclareMathOperator{\Tor}{Tor} \DeclareMathOperator{\QC}{QC}
\DeclareMathOperator{\colim}{colim}

\DeclareMathOperator{\calA}{\mathcal{A}}

\DeclareMathOperator{\calO}{\mathcal{O}}

\DeclareMathOperator{\calB}{\mathcal{B}}

\DeclareMathOperator{\Spec}{Spec}
\DeclareMathOperator{\calF}{\mathcal{F}}
\DeclareMathOperator{\calG}{\mathcal{G}}
\DeclareMathOperator{\Hom}{Hom} 
\DeclareMathOperator{\HH}{H} 
 
\DeclareMathOperator{\calC}{\mathcal{C}}
\DeclareMathOperator{\calI}{\mathcal{I}}

\DeclareMathOperator{\Ind}{Ind} 
\DeclareMathOperator{\calP}{\mathcal{P}} \topmargin=0in

\newtheorem{theorem}{Theorem}[section]
\newtheorem{lemma}[theorem]{Lemma}

\newtheorem{corollary}[theorem]{Corollary}

\theoremstyle{definition}
\newtheorem{definition}[theorem]{Definition}
\newtheorem{example}[theorem]{Example}

\newtheorem{remark}[theorem]{Remark}

\title{Tannaka Duality for Geometric Stacks}
\author{Jacob Lurie}
\begin{document}
\maketitle

\section{Introduction}

Let $X$ and $S$ denote algebraic stacks of finite type over the
field $\C$ of complex numbers, and let $X^{\an}$ and $S^{\an}$
denote their analytifications (which are stacks in the complex
analytic setting). Analytification gives a functor
$$\phi: \Hom_{\C}(S,X) \rightarrow \Hom(S^{\an}, X^{\an}).$$ It is
natural to ask for circumstances under which $\phi$ is an equivalence.

In the case where $X$ and $S$ are projective schemes, a
satisfactory answer was obtained long ago. In this case, both
algebraic and analytic maps may be classified by their graphs,
which are closed in the product $X \times S$. One may then deduce
that any analytic map is algebraic by applying Serre's GAGA
theorem (see \cite{Serre}) to $X \times S$.

If $S$ is a projective scheme and $X$ is the classifying stack of
the algebraic group $\GL_n$, then $\Hom(S,X)$ classifies vector
bundles on $S$. If $S$ is a proper scheme, then any analytic
vector bundle on $S$ is algebraic (again by Serre's GAGA theorem),
and one may again deduce that $\phi$ is an equivalence.

By combining the above methods, one can deduce that $\phi$ is an
equivalence whenever $X$ is given globally as a quotient of a
separated algebraic space by the action of a linear algebraic
group (and $S$ is proper). The main motivation for this paper was
to find a more natural hypothesis on $X$ which forces $\phi$ to be
an equivalence. We will show that this is the case whenever $X$ is
{\it geometric}: that is, when $X$ is quasi-compact and the
diagonal morphism $X \rightarrow X \times X$ is affine. More
precisely, we have the following:

\begin{theorem}\label{decoy}
Let $S$ be a Deligne-Mumford stack which is proper over $\C$, and
let $X$ be a geometric stack of finite type over $\C$. Then the
analytification functor $\phi$ is an equivalence of categories.
\end{theorem}

Our method of proving Theorem \ref{decoy} is perhaps more
interesting than the theorem itself. The basic idea is to show
that if $X$ is a geometric stack, then there exists a Tannakian
characterization for morphisms $f: S \rightarrow X$, in both the
algebraic and analytic categories. More precisely, we will show
that giving a morphism $f$ is equivalent to specifying a
``pullback functor'' $f^{\ast}$ from coherent sheaves on $X$ to
coherent sheaves on $S$. We will then be able to deduce Theorem
\ref{decoy} by applying Serre's GAGA theorem to $S$.

This paper was originally intended to be included in the more
ambitious paper \cite{geometric}, which studies the an analogous
duality theorem in derived algebraic geometry.
However, since the derived setting offers a host of additional
technical difficulties, it seemed worthwhile to write a separate
account in the simpler case considered here.

I would like to thank Brian Conrad for offering several
suggestions and corrections after reading an earlier version of
this paper. I would also like to thank the American Institute of Mathematics for supporting me while this paper was being revised.

\section{Notation}

Throughout this paper, the word {\it topos} shall mean
Grothendieck topos. If $S$ is a topos, then we shall usually refer
to the objects of $S$ as {\it sheaves on $S$}. Similarly we shall
speak of {\it sheaves of groups}, {\it sheaves of rings}, and so
forth, rather than {\it group objects} or {\it ring objects} of
$S$. This terminology is justified by the fact that any category
$S$ is equivalent to the category of representable presheaves on
$S$, and that when $S$ is a topos then the representable
presheaves are precisely those presheaves which are sheaves with
respect to the canonical topology on $S$.

Throughout this note, {\it ring} shall mean {\it commutative ring
with identity}. If $A$ is a ring, then we will write $\Mod_A$ for
the category of $A$-modules. More generally, if $\calA$ is a sheaf
of rings on a topos $S$, then $\Mod_{\calA}$ shall denote the
category of sheaves of $\calA$-modules.

If $\calF$ is an abelian sheaf (on some topos $S$) which admits an
action of a ring of endomorphisms $A$, and $M \in \Mod_A$, then we
shall write $\calF \otimes_A M$ for the {\it external tensor
product} of $\calF$ by $M$ over $A$. In other words, $\calF
\otimes_A M$ is the sheafification (with respect to the canonical
topology) of the presheaf $V \mapsto \calF(V) \otimes_A M$ on $S$.

Throughout this paper, the term {\it algebraic stack} shall mean
Artin stack (not necessarily of finite type) over $\Spec \Z$. If
$X$ is an algebraic stack, we will write $\QC_X$ for the category
of quasi-coherent sheaves on $X$. If $X$ is locally Noetherian,
then we will write $\Coh_X$ for the category of coherent sheaves
on $X$.

\section{Geometric Stacks}

\begin{definition}\label{geom}
An algebraic stack $X$ is {\it geometric} if it is quasi-compact
and the diagonal morphism $$X \rightarrow X \times X = X
\times_{\Spec \Z} X$$ is representable and affine.
\end{definition}

\begin{remark}
The terminology we have just introduced is borrowed from
\cite{toen}, with one modification: we include a hypothesis of
quasi-compactness in our definition of a geometric stack.
\end{remark}

\begin{remark}
Let $X$ be a geometric stack. Since $X$ is quasi-compact, there
exists a smooth surjection $\Spec A \rightarrow X$. Since the
diagonal $X \rightarrow X \times X$ is affine, the fiber product
$\Spec A \times_X \Spec A = \Spec B$ is affine. Moreover, the pair
of objects $(\Spec A, \Spec B)$ are part of a groupoid object in
the category of affine schemes. Algebraically, this means that the
pair of rings $(A,B)$ are endowed with the structure of a {\it
Hopf algebroid}. This Hopf algebroid is commutative (in the sense
that the rings $A$ and $B$ are commutative) and smooth (in the
sense that either of the natural maps $A \rightarrow B$ is
smooth).

Conversely, any commutative, smooth Hopf algebroid gives rise to a
geometric stack in a natural way. It follows that this entire
paper could be written using the language of Hopf algebroids,
rather than algebraic stacks. We will avoid following this course,
since the language of algebraic stacks seems more intuitive and
notationally simpler. We refer the reader to \cite{algebroid} for
a discussion which touches upon some of the ideas of this paper,
written in the language of Hopf algebroids.
\end{remark}

\begin{remark}
In Definition \ref{geom}, we are free to replace the absolute
product $X \times X = X \times_{\Spec \Z} X$ with the fiber
product $X \times_Y X$ for any separated scheme $Y$ which admits a
map from $X$.
\end{remark}

\begin{example}
Any quasi-compact, separated scheme (or algebraic space) is a
geometric stack.
\end{example}

\begin{example}
The classifying stack of any smooth, affine group scheme is a
geometric stack.
\end{example}

\begin{example}
Call a morphism $X \rightarrow S$ of algebraic stacks {\it
relatively geometric} if $X \times_S \Spec A$ is a geometric
stack, for any morphism $\Spec A \rightarrow S$. Then a
composition of relatively geometric morphisms is relatively
geometric. Applying this to the particular case where $S = BG$ is
the classifying stack of an smooth affine group scheme, we deduce
that the quotient of any separated algebraic space by the action
of a smooth affine group scheme is geometric.
\end{example}

The main theme of this paper is that if $X$ is a geometric stack,
then $X$ has ``enough'' quasi-coherent sheaves. As an illustration
of this principle, we prove the following:

\begin{theorem}
Let $X$ be a geometric stack. Then the left-bounded derived
category of quasi-coherent sheaves on $X$ is naturally equivalent
to the full subcategory of the left-bounded derived category of
(smooth-{\et}ale) $\calO_{X}$-modules which have quasi-coherent
cohomologies.
\end{theorem}

\begin{proof}
Let $\QC_X$ denote the abelian category of quasi-coherent sheaves
on $X$. It will suffice to show that $\QC_{X}$ has enough
injective objects, and that if $I \in \QC_{X}$ is injective and $M
\in \QC_{X}$ is arbitrary, then $\Ext^i(M,I) = 0$ for all $i
> 0$. Here the $\Ext$-group is computed in the larger category of all
(smooth-{\et}ale) $\calO_{X}$-modules.

Since $X$ is quasi-compact, we may choose a smooth surjection $p:
U \rightarrow X$, where $U = \Spec A$ is an affine scheme. If $N
\in \QC_{X}$, then we may choose an injection $p^{\ast} N
\rightarrow I$, where $I$ is a quasi-coherent sheaf on $U$
corresponding to an injective $A$-module. Since $p^{\ast}$ is
exact, $p_{\ast} I$ is injective. We claim that the adjoint
morphism $N \rightarrow p_{\ast} I$ is a monomorphism. For this,
it suffices to show that each of the maps $N \rightarrow p_{\ast}
p^{\ast} N$ and $p_{\ast} p^{\ast} N \rightarrow p_{\ast} I$ are
monomorphisms. For the first map, this follows from the fact that
$U$ is a flat covering of $X$. For the second, we note that since
$X$ is geometric, $p$ is an affine morphism so that $p_{\ast}$ is
an exact functor when restricted to quasi-coherent sheaves.

The above argument shows that $\QC_{X}$ has enough injectives, and
that in fact every injective is a direct summand of a
quasi-coherent sheaf having the form $p_{\ast} I$. We remark that since $p$ is a smooth morphism, it induces a geometric morphism from the smooth-{/et}ale topos of $U$ to the smooth-{\et}ale topos of $X$,  and that the quasi-coherent direct-image functor $p_{\ast}$ is the restriction to $\QC_U$ of a functor (also denoted by $p_{\ast}$) defined on {\em all} smooth-{/et}ale sheaves. Now note
that $\Ext^i(M, p_{\ast} I) = \Ext^i( p^{\ast} M, I)$ (since $I$
has vanishing higher direct images under $p$). It now suffices to
show that for any quasi-coherent sheaf $N$ on $U$, we have
$\Ext^i(N,I) =0$ for $i > 0$, where the $\Ext$-group is computed
in the category of smooth-{\et}ale $\calO_U$-modules.

Since $U$ is
affine, there exists a resolution $P_{\bigdot}$ of $N$ such that
each $P_i$ is a direct sum of copies of $\calO_U$. Since
$\Ext^j(\calO_U,I) = \HH^j(U,I) = 0$ for $j > 0$, and since
$\Ext^j(\bigdot, I)$ carries arbitrary direct sums into direct
products, we deduce that $\Ext^i(N,I)$ is the $i$th cohomology
group of the complex $\Hom(P_{\bigdot}, I)$. This cohomology group
vanishes for $i > 0$ since $P_{\bigdot}$ is acyclic in positive
degrees and $I$ is obtained from an injective $A$-module.
\end{proof}

If $X$ is a Noetherian geometric stack, then we can say even more:
$X$ has ``enough'' coherent sheaves. This follows from a
well-known argument, but we include the proof for lack of a
reference:

\begin{lemma}\label{del}
Let $X$ be an algebraic stack which is Noetherian and geometric.
Then $\QC_{X}$ is equivalent to the category of $\Ind$-objects of
the full subcategory $\Coh_X \subseteq \QC_X$.
\end{lemma}

\begin{proof}
We first prove that if $M$ is a coherent sheaf on $X$, then $M$ is
a compact object of $\QC_{X}$. Indeed, suppose that $\{ N_{\alpha}
\}$ is some filtered system of quasi-coherent sheaves on $X$ with
colimit $N$.

Choose a smooth surjection $p: U \rightarrow X$, where $U = \Spec
A$ is affine. Since $p_{\ast}$ and $p^{\ast}$ commute with
filtered colimits, we have a filtered system of short exact
sequences
$$ \{ 0 \rightarrow N_{\alpha} \rightarrow p_{\ast} p^{\ast}
N_{\alpha} \rightarrow p_{\ast} p^{\ast} p_{\ast} p^{\ast}
N_{\alpha} \}$$ having filtered colimit $$ 0 \rightarrow N
\rightarrow p_{\ast} p^{\ast} N \rightarrow p_{\ast} p^{\ast}
p_{\ast} p^{\ast} N.$$ Using these short exact sequences, we see
that in order to prove that $\Hom(M, N) \simeq \colim \{ \Hom(M,
N_{\alpha})\}$, it suffices to prove the analogous result for the
filtered systems $\{ p_{\ast} p^{\ast} N_{\alpha} \}$ and $\{
p_{\ast} p^{\ast} p_{\ast} p^{\ast} N_{\alpha} \}$. In other
words, we may reduce to the case where the filtered system $\{
N_{\alpha} \}$ is the direct image of a filtered system $\{
P_{\alpha} \}$ of quasi-coherent sheaves on $U$; let $P$ be the
colimit of this system. In this case, we have $$\Hom(M, N) =
\Hom(p^{\ast} M, P) = \colim \{ \Hom(p^{\ast} M, P_{\alpha})\} =
\colim \{ \Hom(M, p_{\ast} P_{\alpha}) \}.$$ Here the second
equality follows from the fact that $p^{\ast} M$ is coherent, and
therefore corresponds to a finitely presented $A$-module. We
remark that this last argument also shows that any compact object
of $\QC_{X}$ is coherent.

By formal nonsense, we obtain a fully faithful embedding $\Ind(
\Coh_X) \rightarrow \QC_X$. To complete the proof, it suffices to
show that every quasi-coherent sheaf $M$ on $X$ is a filtered
colimit of coherent subsheaves. For this, we write $p^{\ast} M$ as
a filtered colimit of coherent subsheaves $\{ P_{\alpha} \subseteq
p^{\ast} M\}$ on $U$. Then $M \subseteq p_{\ast} p^{\ast} M =
\bigcup \{ p_{\ast} P_{\alpha} \}$. Set $M_{\alpha} = M \cap
p_{\ast} P_{\alpha}$. Then the natural map $p^{\ast} M_{\alpha}
\rightarrow p^{\ast} M$ factors through $P_{\alpha}$, so it
follows that $M_{\alpha}$ is a coherent subsheaf of $M$. Clearly
$M$ is the union of the filtered family of subobjects $\{
M_{\alpha} \}$.
\end{proof}

\begin{remark}
The proof does not really require that $X$ is geometric; really
all that is needed is that the diagonal morphism $X \rightarrow X
\times X$ is a quasi-compact, quasi-separated relative algebraic
space.
\end{remark}

\section{Maps into Algebraic Stacks}

The main goal of this paper is to prove Theorem \ref{decoy}, which
furnishes a comparison between the categories $\Hom(S,X)$ and
$\Hom(S^{\an}, X^{\an})$. As an intermediate step, we will define
a category $\Hom(S^{\an},X)$. This category will be equivalent to
$\Hom(S^{\an}, X^{\an})$ by construction, and we will be reduced
to comparing $\Hom(S,X)$ with $\Hom(S^{\an}, X)$.

Let $(S, \calO_S)$ be {\em any} ringed topos, and $\calF$ any
covariant functor from commutative rings to groupoids. The
assignment $U \mapsto \calF(\calO_S(U))$ determines a presheaf of
groupoids $\sHom_0(S,\calF)$ on the topos $S$. This presheaf of
groupoids may or may not be a stack on $S$; in either case, there
always exists a stack $\sHom(S,\calF)$ on $S$ which is initial
among stacks equipped with a morphism $$\sHom_0(S,\calF)
\rightarrow \sHom(S,\calF).$$ We will refer to $\sHom(S,\calF)$ as
the {\it stackification} of $\sHom_0(S,\calF)$. The groupoid of
global sections of $\sHom(S,\calF)$ will be denoted by
$\Hom(S,\calF)$. In the case where $\calF$ is represented by an
algebraic stack $X$, we will also write $\sHom(S, X)$ and $\Hom(S,
X)$. We note that this is an abuse of notation, because these
morphism spaces depend on the sheaf of rings $\calO_S$ and not
only on the underlying topos $S$.

\begin{example}
Let $(S, \calO_S)$ be the {\et}ale topos of a Deligne-Mumford
stack, and let $X$ be an arbitrary algebraic stack. Then the
definition of $\Hom(S,X)$ given above agrees with the usual definition.
\end{example}

The definition given above is not of much use unless we have some
means of calculating $\Hom(S,X)$ in terms of a presentation of
$X$. This requires an additional hypothesis on the ringed topos
$(S, \calO_S)$ which we now introduce:

\begin{definition}
A ringed topos $(S, \calO_S)$ is {\it local for the {\et}ale
topology} if it has the following property: for any $E \in S$ and
any finite set of {\et}ale ring homomorphism $\{ \calO_S(E)
\rightarrow R_i \}$, having the property that the induced map
$$\calO_S(E) \rightarrow \prod_i R_i$$ is faithfully flat,
there exist morphisms $E_i \rightarrow E$ in $S$ and
factorizations $\calO_S(E) \rightarrow R_i \rightarrow
\calO_S(E_i)$ having the property that the induced map $$\coprod
E_i \rightarrow E$$ is an epimorphism.
\end{definition}

\begin{remark}
In fact, there exists a canonical choice for $E_i$. For any ring
homomorphism $\calO_S(E) \rightarrow R$, the functor $U \mapsto
\Hom_{\calO_S(E)}(R, \calO_S(U))$ is representable by an object
$U_0 \rightarrow E$. If we imagine that $R$ is presented over
$\calO_S(E)$ by generators and relations, then $U_0$ may be
thought of as the ``sheaf of solutions'' to the corresponding
equations; the requirement of the above definition is that the natural map $E_i \rightarrow E$
be an epimorphism whenever $R$ is {\et}ale over $\calO_S(E)$.
\end{remark}

\begin{remark}
If $S$ has enough points, then $(S, \calO_S)$ is local for the
{\et}ale topology if and only if the stalk $\calO_{S,s}$ at any
point $s$ of $S$ is a strictly Henselian local ring. In
particular, this implies that each stalk $\calO_{S,s}$ is local so
that $(S, \calO_S)$ is a locally ringed topos in the usual sense.

It follows that the {\et}ale topos of a Deligne-Mumford stack is
local for the {\et}ale topology. Similarly, if $(S, \calO_S)$ is
the underlying topos of a complex analytic space, then $(S,
\calO_S)$ is local for the {\et}ale topology.
\end{remark}

\begin{remark}
Suppose that $(S,\calO_S)$ is local for the {\et}ale topology, and
that $X$ is (a functor representable by) a scheme. Then $\Hom(S,X)$ may be
identified with the category of morphisms from $S$ to $X$ in the
$2$-category of locally ringed topoi.
\end{remark}

Let $\calF$ be any groupoid-valued functor defined on commutative
rings, and let $\calF'$ denote the stackification of $\calF$ with
respect to the {\et}ale topology. If $(S, \calO_S)$ is local for
the {\et}ale topology, then the natural map $\sHom(S,\calF)
\rightarrow \sHom(S, \calF')$ is an equivalence of stacks on $S$.

In particular, let us suppose that $X$ is an algebraic stack
equipped with an smooth atlas $p: U \rightarrow X$, so that $(U, U
\times_X U)$ extends naturally to a groupoid object in the
category of algebraic spaces. This groupoid object represents a
functor $\calF$ from rings to groupoids, and $X$ represents the
stackification of the functor $\calF$ with respect to the {\et}ale
topology. If $(S, \calO_S)$ is local for the {\et}ale topology,
then we get $\Hom(S, X) \simeq \Hom(S, \calF)$. It follows that
$\Hom(S,X)$ can be computed in terms of any atlas for $X$. More
concretely, this means that:

\begin{itemize}
\item Locally on $S$, any morphism $f: S \rightarrow X$ factors
through $U$.

\item Given any two morphisms $f,g: S \rightrightarrows U$, any
isomorphism $\alpha: p \circ f \simeq p \circ g$ is induced
locally by a factorization $S \stackrel{h}{\rightarrow} U \times_X
U \rightarrow U \times U$ of $f \times g$.
\end{itemize}

\begin{remark}
If $(S, \calO_S)$ is the underlying topos of a complex analytic
space (or complex-analytic orbifold) and $X$ is any algebraic
stack of finite type over $\C$, then $\Hom_{\C}(S,X) \simeq
\Hom(S,X^{\an})$. To prove this, we note that equality holds when
$X$ is a scheme or algebraic space, essentially by the definition
of the analytification functor. In the general case, both sides
are computed in the same way from a presentation of $X$.
\end{remark}

We conclude this section with a brief discussion of the pullback
functor $f^{\ast}$ determined by a map $f: S \rightarrow X$.
Suppose that $(S, \calO_S)$ is a ringed topos, $X$ an algebraic
stack, and $f: S \rightarrow X$ is any morphism. Locally on $S$,
the morphism $f$ admits a factorization $S \rightarrow \Spec
\Gamma(S,\calO_S) \rightarrow X$, and we may define $f^{\ast}$ as
the composite of the usual pullback functor $\QC_X \rightarrow
\QC_{\Spec \Gamma(S,\calO_S)} = \Mod_{\Gamma(S,\calO_S)}$,
followed by the functor $$M \mapsto \calO_{S} \otimes_{\Gamma(S,
\calO_S)} M$$ from $\Mod_{\Gamma(S, \calO_S)}$ to
$\Mod_{\calO_S}$. This local construction is natural and therefore
makes sense even when $f$ does not factor through $\Spec
\Gamma(S,\calO_S)$. Moreover, the functor $f^{\ast}$ is compatible
with tensor products in the sense that there exist natural
isomorphisms
$$\gamma_{M,N}: f^{\ast}(M \otimes N) \simeq f^{\ast} M \otimes
f^{\ast} N$$
$$\epsilon: f^{\ast} \calO_{X} \simeq \calO_S.$$
The functor $f^{\ast}$ and the coherence data $\{ \gamma_{M,N},
\epsilon \}$ enjoy a number of additional properties which the
next section will place in a more formal context.

\section{Abelian Tensor Categories}

The main step in the proof of Theorem \ref{decoy} is Theorem
\ref{mainevent}, which asserts roughly that a geometric stack $X$
is determined by the category $\QC_{X}$. In order to make a more
precise statement, we must first decide what sort of object
$\QC_{X}$ is. The relevant definitions and the statement of our
main result, Theorem \ref{mainevent}, will be given in this
section.

Recall that a {\it symmetric monoidal category} is a category
$\calC$ equipped with a tensor product bifunctor $\otimes: \calC
\times \calC \rightarrow \calC$ which is coherently unital,
associative, and commutative. This means that there exists an
object $1 \in \calC$ and natural isomorphisms
$$ 1 \otimes M \simeq M \simeq M \otimes 1$$
$$ (M \otimes N) \otimes P \simeq M \otimes (N \otimes P)$$
$$ M \otimes N \simeq N \otimes M.$$
These isomorphisms are required to specify a number of coherence
conditions: for a discussion, we refer the reader to
\cite{symmetric}. These conditions are evidently satisfied in the
cases of relevance to us, and will not play an important role in
this paper.

\begin{remark}
The commutativity and associativity isomorphisms are part of the
data of a symmetric monoidal category. However, we will abuse
notation and simply refer to $(\calC, \otimes)$ or $\calC$ as a
symmetric monoidal category.
\end{remark}

We also recall that a {\it Grothendieck abelian category} is an
abelian category with a generator which satisfies the axiom
(AB$5$) of \cite{grothendieckabelian}: the existence and exactness
of (small) filtered colimits.

\begin{definition}
An {\it abelian tensor category} is a symmetric monoidal category
$(\calC, \otimes)$ with the following properties:

\begin{enumerate}
\item The underlying category $\calC$ is an abelian category.

\item For any fixed object $M \in \calC$, the functor $N \mapsto M
\otimes N$ commutes with finite colimits. Equivalently, the tensor
product operation $\otimes$ is additive and right-exact.
\end{enumerate}

We shall say that $(\calC, \otimes)$ is {\it complete} if $\calC$
is a Grothendieck abelian category and the functor $N \mapsto M
\otimes N$ commutes with {\em all} (small) colimits, for each
fixed object $M \in \calC$.
\end{definition}

\begin{remark}
If $(\calC, \otimes)$ is an abelian tensor category such that the
underlying category $\calC$ is Grothendieck, then $(\calC,
\otimes)$ is complete if and only if for each $M \in \calC$, the
functor $$N \mapsto M \otimes N$$ has a right adjoint
$$ N \mapsto \Hom(M,N).$$ The ``if'' direction is easy
and the reverse implication follows from the adjoint functor
theorem.
\end{remark}

If $(\calC, \otimes)$ is an abelian tensor category and $M \in
\calC$, then we shall say that $M$ is {\it flat} if the functor $N
\mapsto M \otimes N$ is an exact functor. We shall say that
$(\calC, \otimes)$ is {\it tame} if it has the following property:
for any exact sequence $$0 \rightarrow M' \rightarrow M
\rightarrow M'' \rightarrow 0$$ in $\calC$ such that $M''$ is
flat, and any $N \in \calC$, the induced sequence $$0 \rightarrow
M' \otimes N \rightarrow M \otimes N \rightarrow M'' \otimes N
\rightarrow 0$$ is also exact. Any abelian tensor category which
has enough flat objects to set up a theory of flat resolutions is
tame: this follows from vanishing of the group $\Tor_1(M'',N)$. We
will need to work with abelian tensor categories which do not
satisfy the latter condition; however, all of the abelian tensor
categories which we will encounter will be tame.

\begin{lemma}\label{spill}
Let $(\calC, \otimes)$ be a tame abelian tensor category, and let
$$0 \rightarrow M' \rightarrow M \rightarrow M'' \rightarrow 0$$
be an exact sequence in $\calC$. Suppose that $M''$ is flat. Then
$M$ is flat if and only if $M'$ is flat.
\end{lemma}

\begin{proof}
Let $$0 \rightarrow N' \rightarrow N \rightarrow N'' \rightarrow
0$$ be any short exact sequence in $\calC$. A simple diagram chase
shows that $M \otimes N' \rightarrow M \otimes N$ is a
monomorphism if and only if $M' \otimes N' \rightarrow M' \otimes
N$ is a monomorphism.
\end{proof}

An {\it algebra} in $\calC$ is a commutative monoid in $\calC$:
that is, it is an object $A \in \calC$ equipped with a commutative
and associative multiplication $A \otimes A \rightarrow A$ and a
unit $1 \rightarrow A$ (here $1 \in \calC$ denotes the unit for
the tensor product) satisfying the usual identities.

\begin{lemma}\label{fflat}
Let $(\calC, \otimes)$ be a tame abelian tensor category
containing an algebra $A$. The following conditions are
equivalent:

\begin{enumerate}
\item The algebra $A$ is flat, and $A \otimes M = 0$ implies $M
=0$.

\item The unit morphism $u: 1 \rightarrow A$ is a monomorphism,
and the cokernel of $u$ is flat.
\end{enumerate}
\end{lemma}

\begin{proof}
Let us first suppose that $(1)$ is satisfied and prove $(2)$. This
part of the argument will not require the assumption that $\calC$
is tame. In order to prove that $u$ is a monomorphism, it suffices
to prove that $u$ is a monomorphism after tensoring with $A$
(since tensor product with $A$ cannot annihilate the kernel of $u$
unless the kernel of $u$ is zero). But $u \otimes A: A \rightarrow
A \otimes A$ is split by the multiplication $A \otimes A
\rightarrow A$.

A similar argument proves that the cokernel $A'$ of $u$ is flat.
Since $u \otimes A$ is split injective, $A' \otimes A$ is a direct
summand of $A \otimes A$. It follows that $A' \otimes A$ is flat.
Let $$0 \rightarrow N' \rightarrow N \rightarrow N'' \rightarrow
0$$ be any exact sequence in $\calC$. Tensoring with $A'$, we
obtain an exact sequence $$0 \rightarrow K \rightarrow N' \otimes
A' \rightarrow N \otimes A' \rightarrow N'' \otimes A' \rightarrow
0.$$ Since $A' \otimes A$ is flat, we deduce (from the flatness of
$A$) that $K \otimes A = 0$, so that our hypothesis implies that
$K=0$.

Now suppose that $(2)$ is satisfied. Since the cokernel of $u$ is
flat, $A$ is an extension of flat objects of $\calC$ and therefore
flat by Lemma \ref{spill}. Suppose that $A \otimes M = 0$. Since
the cokernel of $u$ is flat, the assumption that $\calC$ is tame
implies that $1 \otimes M \rightarrow A \otimes M$ is a
monomorphism, so that $M \simeq 1 \otimes M \simeq 0$.
\end{proof}

\begin{definition}
Let $(\calC, \otimes)$ be a tame abelian tensor category. An
algebra $A \in \calC$ is {\it faithfully flat} if the equivalent
conditions of Lemma \ref{fflat} are satisfied.
\end{definition}

We now give some examples of abelian tensor categories.

\begin{example}
Let $(S, \calO_S)$ be a ringed topos. Then the usual tensor
product operation endows the category $\Mod_{\calO_S}$ with the
structure of a complete abelian tensor category. Moreover,
$\Mod_{\calO_S}$ is tame. This follows from the fact that
$\calO_S$ has enough flat sheaves to set up a good theory of
$\Tor$-functors.

If $S$ has enough points, then a sheaf of $\calO_S$-modules
($\calO_S$-algebras) is flat (faithfully flat) if and only its
stalk at every point $s \in S$ is flat (faithfully flat) as an
$\calO_{S,s}$-module (algebra).
\end{example}

\begin{example}\label{goof}
Let $X$ be an algebraic stack. Then the category $\QC_X$, equipped
with its usual tensor structure, is a complete abelian tensor
category. However, we must distinguish between two potentially
different notions of flatness. We will call an object $M \in
\QC_X$ {\it globally flat} if it is flat in the sense defined
above: that is, $N \mapsto M \otimes N$ is an exact functor from
$\QC_X$ to itself. We shall call an object $M \in \QC_X$ {\it
locally flat} if it is flat in the usual algebro-geometric sense:
that is, for morphism $f: \Spec A \rightarrow X$, the $A$-module
$\Gamma(\Spec A, f^{\ast} M)$ is flat. It is easy to see that any
locally flat module is globally flat, but the converse is unclear.

If $X$ is geometric, then any globally flat object $M \in \QC_X$
is locally flat. To prove this, let us choose a smooth surjection
$p: U \rightarrow X$, where $U = \Spec A$ is affine. To show that
$M$ is locally flat, we need to show that $p^{\ast} M$ is flat as
a quasi-coherent sheaf on $U$. In other words, we need to show
that the functor $N \mapsto p^{\ast} M \otimes N$ is an exact
functor from $\QC_U$ to itself. It suffices to prove the exactness
after composing with the pullback functor $\QC_{U} \rightarrow
\QC_{U \times_X U}$. Using the appropriate base-change formula,
this is equivalent to the assertion that the functor
$$N \mapsto p^{\ast} p_{\ast} ( p^{\ast} M \otimes N)$$ is exact.
Making use of the natural push-pull isomorphism $ p_{\ast} (
p^{\ast} M \otimes N) \simeq M \otimes p_{\ast} N$, we are reduced
to proving the exactness of the functor
$$N \mapsto p^{\ast} ( M \otimes p_{\ast} N).$$ This is clear,
since the functor is a composite of the exact functors $p_{\ast}$,
$p^{\ast}$, and $M \otimes \bigdot$.

From the equivalence of local and global flatness, we may deduce
that $\QC_{X}$ is tame whenever $X$ is a geometric stack.
\end{example}

Our next goal is to describe the appropriate notion of functor
between abelian tensor categories.

\begin{definition}
Let $(\calC, \otimes)$ and $(\calC', \otimes')$ be abelian tensor
categories. An {\it additive tensor functor} $F^{\ast}$ from
$\calC$ to $\calC'$ is a symmetric monoidal functor (that is, a
functor which is compatible with the symmetric monoidal structures
on $\calC$ and $\calC'$ up to natural isomorphism; see
\cite{symmetric} for a discussion) which commutes with finite
colimits (this latter condition is equivalent to the condition
that $F^{\ast}$ be additive and right-exact).

If $(\calC, \otimes)$ and $(\calC', \otimes')$ are complete, then
we shall say that $F^{\ast}$ is {\it continuous} if it commutes
with all colimits.

We shall say that $F^{\ast}$ is {\it tame} if it possesses the
following additional properties:

\begin{itemize}
\item If $M \in \calC$ is flat, then $F^{\ast} M \in \calC'$ is
flat.

\item If $$0 \rightarrow M' \rightarrow M \rightarrow M''
\rightarrow 0$$ is a short exact sequence in $\calC$ and $M''$ is
flat, then the induced sequence $$0 \rightarrow F^{\ast} M'
\rightarrow F^{\ast} M \rightarrow F^{\ast} M'' \rightarrow 0$$ is
exact in $\calC'$.
\end{itemize}
\end{definition}

\begin{remark}\label{reff}
Let $F^{\ast}$ be an additive tensor functor between abelian
tensor categories $\calC$ and $\calC'$. Since $F^{\ast}$ is a
symmetric monoidal functor, it carries algebra objects in $\calC$
to algebra objects in $\calC'$. If, in addition, $F^{\ast}$ is
tame, then it carries faithfully flat algebras in $\calC$ to
faithfully flat algebras in $\calC'$: this is clear from the
second characterization given in Lemma \ref{fflat}.
\end{remark}

If $\calC$ and $\calC'$ are complete, tame, abelian tensor
categories, then we shall let $\Hom_{\otimes}(\calC, \calC')$
denote the groupoid of continuous, tame, additive tensor functors
from $\calC$ to $\calC'$ (where the morphisms are given by
isomorphisms of symmetric monoidal functors). It is a full
subcategory of the groupoid of all monoidal functors from $\calC$
to $\calC'$. We remark that the notation is slightly abusive: the
category $\Hom_{\otimes}(\calC, \calC')$ depends on the symmetric
monoidal structures on $\calC$ and $\calC'$, and not only on the
underlying categories.

If $(S, \calO_S)$ is a ringed topos, $X$ an algebraic stack, and
$f: S \rightarrow X$ is any morphism, then in the last section we
constructed an associated pullback functor $f^{\ast}: \QC_X
\rightarrow \Mod_{\calO_S}$. From the local description of
$f^{\ast}$, it is easy to see that $f^{\ast}$ is a continuous,
tame, additive tensor functor from $\QC_{X}$ to $\Mod_{\calO_S}$.
We are now prepared to state the main result of this paper:

\begin{theorem}\label{mainevent}
Suppose that $(S, \calO_S)$ is a ringed topos which is local for
the {\et}ale topology and that $X$ is a geometric stack. Then the
functor $$f \mapsto f^{\ast}$$ induces an equivalence of
categories
$$ T: \Hom(S,X) \rightarrow \Hom_{\otimes}(\QC_X, \Mod_{\calO_S}).$$
\end{theorem}

The proof of Theorem \ref{mainevent} will occupy the next four
sections of this paper.

\begin{remark}
Let $X$ and $S$ be arbitrary algebraic stacks, and define
$\Hom'_{\otimes}(\QC_X, \QC_S) \subseteq
\Hom_{\otimes}(\QC_X,\QC_S)$ to be the full subcategory consisting
of tensor functors which carry flat objects of $\QC_X$ to {\em
locally flat} objects of $\QC_S$ (that is, objects of $\QC_S$
which are flat according the usual definition). In particular, if
every flat object of $\QC_{S}$ is locally flat (for example, if
$S$ is geometric), then $\Hom'_{\otimes}(\QC_X, \QC_S) =
\Hom_{\otimes}(\QC_X, \QC_S)$.

We note that $\Hom(S,X)$ and $\Hom'_{\otimes}(\QC_X, \QC_S)$ are
both stacks with respect to the smooth topology on $S$.
Consequently, to prove that $\Hom(S,X) \simeq
\Hom'_{\otimes}(\QC_X, \Mod_{\calO_S})$ we may work locally on $S$
and thereby reduce to the case where $S$ is affine scheme. In this
case, the result follows from Theorem \ref{mainevent}, at least
when $X$ is geometric. Consequently, Theorem \ref{mainevent} and
Example \ref{goof} imply that the functor
$$ X \mapsto \QC_X$$
is a fully faithful embedding of the $2$-category of geometric
stacks into the $2$-category of tame, complete abelian tensor
categories.

Unfortunately, it seems very difficult to say anything about the
essential image of this functor: that is, to address the question
of when an abelian tensor category arises as the category of
quasi-coherent sheaves on a geometric stack.
\end{remark}

\section{The Proof that $T$ is Faithful}

In this section we give the argument for the first and easiest
step in the proof of Theorem \ref{mainevent}: showing that $T$ is
faithful. Since $\Hom(S,X)$ is a groupoid, this reduces to the
following assertion: if $F: S \rightarrow X$ is any morphism, and
$\alpha$ any automorphism of $F$ such that the natural
transformation $T\alpha: F^{\ast} \rightarrow F^{\ast}$ is the
identity, then $\alpha$ is the identity.

Let $(S, \calO_S)$ be any ringed topos and $p: U \rightarrow X$ a
morphism of algebraic stacks, and let $F: S \rightarrow X$ be any
morphism. We note that the $\calO_X$-algebra morphism $\calO_X
\rightarrow p_{\ast} \calO_U$ acquires a canonical section after
pullback to $U$. We deduce the existence of a natural map $\theta$
from the set of factorizations $\{ f: S \rightarrow U | p \circ f
= F \}$ to the set of sections of the algebra homomorphism
$\calO_S \rightarrow F^{\ast} p_{\ast} \calO_U$. The crucial
observation is the following:

\begin{lemma}
If $(S, \calO_S)$ is local for the {\et}ale topology and $p$ is
affine, then $\theta$ is bijective.
\end{lemma}

\begin{proof}
The assertion is local on $S$. We may therefore suppose that $F$
factors through some smooth morphism $V \rightarrow X$, where $V$
is an affine scheme. Replacing $X$ by $V$ and $U$ by $U \times_X
V$ (and noting that the formation of $p_{\ast}$ is compatible with
the flat base change $V \rightarrow X$), we may reduce to the case
in which $X$ and $U$ are affine schemes. In this situation, the
result is obvious.
\end{proof}

Let us now return to the setting of Theorem \ref{mainevent}. Since
$X$ is quasi-compact, there exists a smooth surjection $p: U
\rightarrow X$, where $U$ is an affine scheme. Since the diagonal
of $X$ is affine, $p$ is an affine morphism. Let $\calA = F^{\ast}
p_{\ast} \calO_U$. The condition that $\alpha$ be the identity is
local on $S$; since $p$ is surjective, we may suppose the
existence of a factorization $S \stackrel{f}{\rightarrow} U
\stackrel{p}{\rightarrow} X$ for $F$. Let $\overline{f} = \theta
f: \calA \rightarrow \calO_S$ be the morphism of sheaves of
algebras classifying $f$.

The morphism $\alpha$ induces a factorization $S \rightarrow U
\times_{X} U \rightarrow X$, which is classified by the
$\calO_S$-algebra map $$\calA \otimes_{\calO_S} \calA \stackrel{1
\otimes T\alpha(\calA)} \rightarrow \calA \otimes_{\calO_S} \calA
\stackrel{\overline{f} \otimes \overline{f}}{\rightarrow}
\calO_S.$$ If $T\alpha$ is the identity, then this
$\calO_S$-algebra map coincides with the map $\theta( \Delta \circ
f)$ classifying the composition $$S \stackrel{f}{\rightarrow} U
\stackrel{\Delta}{\rightarrow} U \times_X U.$$ Since $\theta$ is
injective, we deduce that $\alpha$ is the identity.

\section{The Proof that $T$ is Full}

Our next goal is to prove that the functor $T$ is full.
Concretely, this means that given any pair of morphisms $F,G: S
\rightrightarrows X$ and any isomorphism $\beta: F^{\ast}
\rightarrow G^{\ast}$, there exists an isomorphism $\alpha: F
\rightarrow G$ with $T\alpha = \beta$.

Since we have already shown that $\alpha$ is uniquely determined,
it suffices to construct $\alpha$ locally on $S$. We may therefore
suppose that $F$ factors as $S \stackrel{f}{\rightarrow} U
\stackrel{p}{\rightarrow} X$, where $U = \Spec A$ is affine and
$p$ is a smooth surjection. Similarly, we may suppose that $G$
factors as $S \stackrel{g}{\rightarrow} V
\stackrel{q}{\rightarrow} X$, where $V = \Spec B$ is affine and
$q$ is a smooth surjection. Of course, we could take $U = V$ and
$p=q$, but this would lead to unnecessary confusion.

Let $\calA = F^{\ast} p_{\ast} \calO_U$ and $\calB = G^{\ast}
q_{\ast} \calO_V$. Then $\calA$ and $\calB$ are sheaves of
$\calO_S$-algebras, and $$\calA \otimes \calB
\stackrel{\beta}{\simeq} G^{\ast} p_{\ast} \calO_U \otimes
G^{\ast} q_{\ast} \calO_V \simeq G^{\ast} r_{\ast} \calO_W,$$
where $r: W = U \times_X V \rightarrow X$ is the natural
projection. The sections $f$ and $g$ induce morphisms $\calA
\rightarrow \calO_S$, $\calB \rightarrow \calO_S$ of sheaves of
algebras. Tensoring them together, we obtain a morphism $\eta:
G^{\ast} r_{\ast} \calO_W \rightarrow \calO_S$ which classifies a
morphism $h: S \rightarrow W$. It is clear from the construction
that $h$ induces an isomorphism $\alpha: F \rightarrow G$.

To complete the proof, it will suffice to show that $T\alpha =
\beta$. In other words, we must show that for any $M \in \QC_X$,
the induced maps $T\alpha(M), \beta(M): F^{\ast} M \rightarrow
G^{\ast} M$ coincide. Since $F$ factors through $r: W \rightarrow
X$, the sheaf $F^{\ast} M$ is a direct factor of $F^{\ast}
r_{\ast} r^{\ast} M$. Therefore we may suppose that $M = r_{\ast}
N$ for some $N \in \QC_W$. Choosing a surjection $N' \rightarrow
N$ with $N'$ free, we may reduce to the case where $N$ is free
(since $F^{\ast}$ is right exact). Using the fact that $F^{\ast}$
commutes with direct sums, we may reduce to the case where $N =
\calO_W$. In this case, $F^{\ast} M \simeq \calA \otimes_A C$,
$G^{\ast} M \simeq \calB \otimes_B C$. We now observe that both
$T\alpha$ and $\beta$ are implemented by the isomorphism $\calA
\otimes_A C \simeq \calA \otimes_{\calO_S} \calB \simeq \calB
\otimes_B C$.

\section{Interlude}

Before we can complete the proof of Theorem \ref{mainevent}, we
need to introduce some additional terminology and establish some
lemmas.

\begin{definition}
Let $S$ be a topos, $\calA$ a sheaf of rings on $S$. A sheaf of
$\calA$-modules $\calF$ is {\it locally finitely presented} if $S$
admits a covering by objects $\{ U_{\alpha} \}$ such that each
$\calF|U_{\alpha}$ is isomorphic to $\calA \otimes_{
\calA(U_{\alpha}) } F_{\alpha}$ for some finitely presented
$\calA(U_{\alpha})$-module $F_{\alpha}$. We will say that $\calF$
is {\it locally projective} if each $F_{\alpha}$ may be chosen to
be a (finitely generated) projective module over
$\calA(U_{\alpha})$.

Similarly, if $\calB$ is a
sheaf of $\calA$-algebras on $S$, then we shall say that $\calB$ is
{\it smooth} over $\calA$ if
$S$ admits a
covering by objects $\{ U_{\alpha} \}$ such that $\calB|U_{\alpha}
\simeq \calA|U_{\alpha} \otimes_{ \calA(U_{\alpha}) } R_{\alpha}$,
where each $R_{\alpha}$ is a smooth
$\calA(U_{\alpha})$-algebra.
\end{definition}

\begin{lemma}\label{rem}
Suppose that $(S,\calO_S)$ is a ringed topos which is local for
the {\et}ale topology. Let $\calA$ be a sheaf of
$\calO_S$-algebras which is smooth and faithfully flat over
$\calO_S$. Then, locally on $S$, there exists a section $\calA
\rightarrow \calO_S$.
\end{lemma}

\begin{proof}
Since the assertion is local on $S$, we may suppose that $\calA
\simeq \calO_S \otimes_{\Gamma(S,\calO_S)} A$, where $A$ is a
smooth $\Gamma(S,\calO_S)$-algebra. Consequently, the map $\Spec A
\rightarrow \Spec \Gamma(S, \calO_S)$ is open. Thus, there exist
finitely many global sections $\{ s_1, \ldots, s_n \}$ of
$\calO_S$ such that $A[\frac{1}{s_i}]$ is faithfully flat over
$\Gamma(S, \calO_S)[\frac{1}{s_i}]$, and $$A \rightarrow
A[\frac{1}{s_1}] \times \ldots \times A[\frac{1}{s_n}]$$ is
faithfully flat. Let $\calI$ denote the ideal sheaf of $\calO_S$
generated by $\{ s_1, \ldots, s_n \}$. Then $\calO_S / \calI$ is
annihilated by tensor product with $A[ \frac{1}{s_i}]$ for each
$i$, and therefore annihilated by tensor product with $A$. Since
$\calA$ is faithfully flat, we deduce that $\calO_S / \calI = 0$.
Thus, the global sections $\{ s_1, \ldots, s_n \}$ generate the
unit ideal sheaf. Shrinking $S$ further, we may suppose that the
$s_i$ generate the unit ideal in $\Gamma(S,\calO_S)$. This implies
that $\Spec A \rightarrow \Spec \Gamma(S, \calO_S)$ is surjective.
Since $A$ is smooth over $\Gamma(S, \calO_S)$, we deduce the
existence of a section $A \rightarrow R$, where $R$ is {\et}ale
and faithfully flat over $\Gamma(S, \calO_S)$. Since $\calO_S$ is
local for the {\et}ale topology, we may (after shrinking $S$)
assume the existence of a section $R \rightarrow \Gamma(S,
\calO_S)$. The composite homomorphism $A \rightarrow \Gamma(S,
\calO_S)$ induces the desired section $\calA \rightarrow \calO_S$.
\end{proof}

\begin{lemma}\label{bbq}
Let $\calA$ be a sheaf of rings on a topos $S$, and let $\calF$ be
a sheaf of $\calA$-modules. Then $\calF$ is locally projective if
and only if it is locally finitely presented and, for each $U \in
S$, the functor $\sHom(\calF|U, \bigdot)$ is an exact functor from
$\Mod_{\calA|U}$ to itself.
\end{lemma}

\begin{proof}
The ``only if'' direction is obvious. For the converse, we may
locally choose a surjection $\calA^n \rightarrow \calF$. The
hypothesis implies that the identity map $\calF \rightarrow \calF$
admits a lifting $\calF \rightarrow \calA^n$, at least locally on
$S$, so we may write $\calA^n \simeq \calF \oplus \calF'$.
Consequently, $P=\Gamma(S,\calF)$ is a direct summand of
$\Gamma(S, \calA^n)$ and is therefore a finitely generated,
projective $\Gamma(S,\calA)$-module; let $P' = \Gamma(S, \calF')$
denote the complementary factor. Let $\calP = \calA
\otimes_{\Gamma(S, \calA)} P$ and $\calP' = \calA
\otimes_{\Gamma(S, \calA)} P$.

The isomorphisms $P \rightarrow \Gamma(S,\calF)$ and $P'
\rightarrow \Gamma(S, \calF')$ induce maps $\alpha: \calP
\rightarrow \calF$ and $\alpha': \calP' \rightarrow \calF'$. The
direct sum $\alpha \oplus \alpha'$ is the isomorphism $$\calP
\oplus \calP' \simeq \calA \otimes_{\Gamma(S, \calA)} (P \oplus
P') \simeq \calA^n.$$ It follows that $\alpha$ and $\alpha'$ are
both isomorphisms, so that $\calF \simeq \calP$ is locally
projective.
\end{proof}

\begin{lemma}\label{bb}
Let $S$ be a topos, $\calA$ a sheaf of rings on $S$, $\calB$ a
faithfully flat sheaf of $\calA$-algebras, and $\calF$ a locally
finitely presented sheaf of $\calA$-modules. Then $\calF$ is
locally projective (as a sheaf of $\calA$-modules) if and only if
$\calF \otimes_{\calA} \calB$ is locally projective (as a sheaf of
$\calB$-modules).
\end{lemma}

\begin{proof}
The ``only if'' direction is obvious. For the converse, suppose
that $\calF \otimes_{\calA} \calB$ is locally projective. Since
$\calF$ is locally finitely presented, the assumption that $\calB$
is flat over $\calA$ implies that $$\sHom_{\Mod_{\calA}}(\calF,
\calG) \otimes_{\calA} \calB \simeq \sHom_{\Mod_{\calB}}(\calF
\otimes_{\calA} \calB, \calG \otimes_{\calA} \calB).$$ The
exactness of the latter functor implies the exactness of the
former, since $\calB$ is faithfully flat over $\calA$. Now we
simply apply Lemma \ref{bbq}.
\end{proof}

\begin{lemma}\label{treb}
Let $(S, \calO_S)$ be a ringed topos, let $f: \calB' \rightarrow
\calB$ be a map of $\calO_S$-algebras, and let $\calA$ be a
faithfully flat $\calO_S$-algebra. Suppose that the map $f
\otimes_{\calO_S} \calA: \calB' \otimes_{\calO_S} \calA \simeq
\calB \otimes_{\calO_S} \calA$ extends to an isomorphism
$$ \calB' \otimes_{\calO_S} \calA \simeq (\calB \otimes_{\calO_S}
\calA) \times \calC$$ for some $\calA$-algebra $\calC$. Then $f$
extends to an isomorphism $$\calB' \simeq \calB \times \calC_0$$
of $\calO_S$-algebras.
\end{lemma}

\begin{proof}
Since $f \otimes_{\calO_S} \calA$ is surjective, $f$ is
surjective. Let $\calI \subseteq \calB'$ denote the kernel of $f$.
The hypothesis implies that $\calI \otimes_{\calO_S} \calA$ is
generated by an idempotent. It follows that for any
$\calA$-algebra $\calA'$, the ideal sheaf $\calI \otimes_{\calO_S}
\calA'$ is generated by a idempotent section
$$ e_{\calA'} \in \Gamma(S, \calI \otimes_{\calO_S} \calA').$$
We note that the idempotent $e_{\calA'}$ is uniquely determined
and independent of the $\calA$-algebra structure on $\calA'$. In
particular, we deduce that the image of $e_{\calA}$ under the two
natural maps $$\Gamma(S,\calI \otimes_{\calO_S} \calA)
\rightrightarrows \Gamma(S,\calI \otimes_{\calO_S} \calA
\otimes_{\calO_S} \calA)$$ coincide. Since $\calA$ is faithfully
flat, we deduce that $e_{\calA}$ belongs to the image of the
injection $\Gamma(S,\calI) \rightarrow \Gamma(S,\calI
\otimes_{\calO_S} \calA)$. Let us denote its preimage by $e \in
\Gamma(S,\calI)$. One deduces readily that $e$ is an idempotent
which generates $\calI$, which gives rise to the desired product
decomposition for $\calB'$.
\end{proof}

\begin{lemma}\label{trub}
Let $(S, \calO_S)$ be a ringed topos, and let $\calA$ and $\calB$
be sheaves of $\calO_S$-algebras. Suppose that $\calA$ is
faithfully flat over $\calO_{S}$ and that $\calA \otimes_{\calO_S}
\calB$ is smooth over $\calA$. Then, locally on $S$, we may find
finitely many global sections $\{ b_1, \ldots, b_n \}$ of $\calB$
which generate the unit ideal of $\Gamma(S,\calB)$, with the
property that each $\calB[\frac{1}{b_i}]$ is smooth over
$\calO_S$.
\end{lemma}

\begin{proof}
Let $\Omega$ denote the sheaf of relative differentials of $\calB$
over $\calO_S$. In other words, $\Omega$ is the sheafification of
the presheaf $U \mapsto \Omega_{\calB(U)/\calO_S(U)}$, where
$\Omega_{R/R'}$ denotes the module of relative differentials of
$R$ over $R'$. In particular, if $\calB = \calO_S
\otimes_{\Gamma(S, \calO_S)} R$, then $\Omega = \calO_S
\otimes_{\Gamma(S, \calO_S)} \Omega_{R/\Gamma(S, \calO_S)}$.

The sheaf $\Omega \otimes_{\calB} (\calA \otimes_{\calO_S} \calB)$
is isomorphic to the sheaf of relative differentials of $\calA
\otimes_{\calO_S} \calB$ over $\calA$, and therefore locally
projective. Consequently, $\Omega$ is locally projective as a
$\calB$-module by Lemma \ref{bb}.

We may suppose that $\calA \otimes_{\calO_S}
\calB$ is generated over $\calA$ by finitely many global section
$\{x_1, \ldots, x_n\}$ of $\calB$. The corresponding map
$\calO_S[x_1, \ldots, x_n] \rightarrow \calB$ becomes surjective
after tensoring with $\calA$. Since $\calA$ is faithfully flat, we
deduce that $\calO_{S}[x_1, \ldots, x_n] \rightarrow \calB$ is
surjective. Thus, the induced map on differentials
$\calB^n \rightarrow \Omega$ is surjective. Let $P$ be the kernel of this
map, so that $P$ is a locally projective $\calB$-module.

Locally on $S$, we may find sections $\{ b_1, \ldots, b_m\}$ of
$\calB$ such that $P[\frac{1}{b_i}]$ is a finitely generated, free
$\calB[\frac{1}{b_i}]$-module. Replacing $\calB$ by
$\calB[\frac{1}{b_i}]$, we may reduce to the case where $P$ is
free of rank $(n-k)$.

Let $\calI$ denote the kernel of the induced map $\calO_S[x_1,
\ldots, x_n] \rightarrow \calB$. Then $\calI/\calI^2 \simeq P$.
Locally on $S$, we may find sections $\{ f_1, \ldots, f_k \}$ of
$\calI$ which freely generate $P$ over $\calB$. Localizing $S$
further, we may suppose that $f_1, \ldots, f_k \in
\Gamma(S,\calO_S)[x_1, \ldots, x_n]$. Let $a_1, \ldots, a_m \in
\Gamma(S,\calO_S)[x_1, \ldots, x_n]$ denote the determinants of $k
\times k$ minors of the Jacobian matrix for $(f_1, \ldots, f_k)$.
The assumption that the differentials $\{ df_i \}$ generate a
summand of $\calB^n$ implies that $a_1, \ldots, a_m$ generate the
unit ideal sheaf in $\calB$. Consequently, working locally on $S$
and replacing $\calB$ by $\calB[ \frac{1}{a_i}]$, we may suppose
that some $a_i$ is invertible in $\calB$.

Let $\calB' = \calO_{S}[x_1, \ldots, x_n][\frac{1}{a_i}]/(f_1,
\ldots, f_k)$. Then, by construction, $\calB'$ is smooth and we
have a surjection $p: \calB' \rightarrow \calB$. To prove that
$\calB$ is smooth, it suffices to show that $p$ can be extended to
an isomorphism $\calB' \simeq \calB \times \calC_0$. By Lemma
\ref{treb}, it will suffice to show that $f \otimes_{\calO_S}
\calA$ can be extended to an isomorphism $\calB' \otimes_{\calO_S}
\calA \simeq (\calB \otimes_{\calO_S} \calA) \times \calC$. The
latter assumption is local on $S$, so we may suppose that $\calB'
\otimes_{\calO_S} \calA \simeq \calA \otimes_{\Gamma(S,\calA)} R'$
and $\calB \otimes_{\calO_S} \calA \simeq \calA
\otimes_{\Gamma(S,\calA)} R$, where $R$ and $R'$ are smooth
$\Gamma(S,\calA)$-algebras of relative dimension $(n-k)$ over
$\Gamma(S,\calA)$, and that $f \otimes_{\calO_S} \calA$ is induced
by a surjection $g: R' \rightarrow R$.

We note that $g$ induces a closed immersion $\Spec R \rightarrow
\Spec R'$ of smooth $\Spec \Gamma(S,\calA)$-schemes having the
same relative dimension over $\Spec \Gamma(S,\calA)$. It follows
that this closed immersion is also an open immersion, so that $g$
extends to an isomorphism $R' \simeq R \times R''$, which
evidently gives rise to the desired factorization of $\calB
\otimes_{\calO_S} \calA$.
\end{proof}

\begin{remark}
It seems likely that it is possible to take $\{ b_1, \ldots, b_n \} =
\{1\}$ in the statement of Lemma \ref{trub}. However, we were unable to
prove this without the assumption that $S$ has enough points. Lemma
\ref{trub} will be sufficient for our application.
\end{remark}

\section{The Proof that $T$ is Essentially Surjective}

In this section, we will complete the proof of Theorem
\ref{mainevent} by showing that the functor $T$ is essentially
surjective. In other words, we must show that if $F^{\ast}: \QC_X
\rightarrow \Mod_{\calO_S}$ is a tame, continuous, additive tensor
functor, then $F^{\ast}$ is the pullback functor associated to
some morphism $F: S \rightarrow X$. Since we have already shown
that $T$ is fully faithful, the morphism $F$ is uniquely
determined; it therefore suffices to construct $F$ locally on $S$.

Let $p: U \rightarrow X$ be a smooth surjection, where $U = \Spec
A$ is affine. Let $\calA = F^{\ast} p_{\ast} \calO_U$. Since
$F^{\ast}$ is an additive tensor functor, $\calA$ is a sheaf of
$\calO_S$-algebras. Moreover, the isomorphism $A \simeq \Gamma(U,
\calO_U)$ induces a morphism $A \rightarrow \Gamma(S, \calA)$. Let
$U \times_X U \simeq \Spec B$ so that $p_{\ast} \calO_U
\otimes_{\calO_X} p_{\ast} \calO_U \simeq p_{\ast} \calO_U
\otimes_A B$. Since $F^{\ast}$ commutes with all colimits, it
commutes with external tensor products, so that $\calA
\otimes_{\calO_S} \calA \simeq \calA \otimes_A B$. Since $B$ is
smooth over $A$, we deduce that $\calA \otimes_{\calO_S} \calA$ is
smooth over $\calA$. Since $F^{\ast}$ is a tame functor, $\calA$
is faithfully flat over $\calO_S$ (see Remark \ref{reff}).
Applying Lemma \ref{trub}, we deduce (possibly after shrinking
$S$) the existence of finitely many global sections $\{ a_1,
\ldots, a_n\}$ of $\calA$, which generate the unit ideal of
$\Gamma(S,\calA)$, such that $\calA[\frac{1}{a_i}]$ is smooth. Let
$\calA'$ denote the product of the $\calO_S$-algebras
$\calA[\frac{1}{a_i}]$. Then $\calA'$ is smooth and faithfully
flat over $\calO_S$, so there exists a section $\calA' \rightarrow
\calO_S$. Composing with the natural map $\calA \rightarrow
\calA'$, we deduce the existence of a section $s: \calA
\rightarrow \calO_S$.

The composite map $A \rightarrow \Gamma(S,\calA)
\stackrel{s}{\rightarrow} \Gamma(S,\calO_S)$ induces a morphism
$g: S \rightarrow \Spec A = U$. We claim that the composition $G =
p \circ g: S \rightarrow X$ has the desired properties. To prove
this, we must exhibit an isomorphism $G^{\ast} \simeq F^{\ast}$ of
additive tensor functors.

Let $M$ be a quasi-coherent sheaf on $X$. Then $p_{\ast} p^{\ast}
M \simeq (p_{\ast} \calO_U) \otimes_A \Gamma(U,M)$, so that
$F^{\ast} p_{\ast} p^{\ast} M \simeq \calA \otimes_A \Gamma(U,M)$.
Composing with the adjunction morphism $M \rightarrow p_{\ast}
p^{\ast} M$ and using the section $s: \calA \rightarrow \calO_S$,
we deduce the existence of a natural transformation $\beta_M:
F^{\ast} M \rightarrow \calO_S \otimes_A M = G^{\ast} M$. It is
easy to see that this is a map of additive tensor functors. To
complete the proof, it suffices to show that $\beta_M$ is an
isomorphism for every $M$.

Since $\calA$ is faithfully flat, to show that $\beta_M$ is an
isomorphism it suffices to show that $\beta_M \otimes_{\calO_S}
\calA = \beta_{M \otimes p_{\ast} \calO_U}$ is an isomorphism. In
particular, we may suppose that $M = p_{\ast} N$ for some
$\calO_U$-module $N$. Since both $F^{\ast}$ and $G^{\ast}$ are
right exact and commute with all direct sums, we may reduce to the
case where $N = \calO_U$. In this case one can easily compute that
$F^{\ast} M \simeq \calA \simeq G^{\ast} M$ and that $\beta_M$
corresponds to the identity map.

\section{The Proof of Theorem \ref{decoy}}

The goal of this section is to show that Theorem \ref{mainevent}
implies Theorem \ref{decoy}. The main difficulty that needs to be
overcome is that Theorem \ref{mainevent} is concerned with
categories of quasi-coherent sheaves, which are rather unwieldy.
We therefore specialize to the case where $X$ is a {\it
Noetherian} geometric stack. In this case, we have a well-behaved
subcategory of coherent sheaves $\Coh_X \subseteq \QC_X$. We let
$\Coh_S \subseteq \Mod_{\calO_S}$ denote the category of locally
finitely presented $\calO_S$-modules. If $\calO_S$ is a coherent
sheaf of rings in the usual sense, then this agrees with the usual
notion of a coherent $\calO_S$-module. If $\calO_S$ is not
coherent, then $\Coh_S$ need not be an abelian subcategory of
$\Mod_{\calO_S}$, but this does not impact any of the statements
which follow.

By Lemma \ref{del}, the inclusion $\Coh_X \subseteq \QC_X$ induces
an equivalence of categories $\Ind(\Coh_X) \simeq \QC_X$. It
follows that the category of continuous, additive tensor functors
from $\QC_X$ to $\Mod_{\calO_S}$ is equivalent to the category
$\calC$ of additive tensor functors from $\Coh_X$ to
$\Mod_{\calO_S}$. By Theorem \ref{mainevent}, $\Hom(S,X)$ is
equivalent to the full subcategory of $\calC$, consisting of those
functors $\Coh_X \rightarrow \Mod_{\calO_S}$ which admit
continuous, tame extensions to $\QC_X$. Since Theorem
\ref{mainevent} implies that any such functor has the form
$f^{\ast}$ for some map $S \rightarrow X$, it must carry coherent
$\calO_X$-modules to coherent $\calO_S$-modules (this is immediate
from the construction of $f^{\ast}$). Consequently, we may deduce
the following ``coherent'' version of Theorem \ref{mainevent}:

\begin{corollary}\label{corevent}
Suppose that $(S, \calO_S)$ is ringed topos which is local for the
{\et}ale topology, such that $\calO_S$ is coherent. Let $X$ be a
Noetherian geometric stack, and let $\Hom(\Coh_X, \Coh_S)$ denote
the groupoid of additive tensor functors from $\Coh_X$ to
$\Coh_S$. Then the natural functor
$$ \Hom(S,X) \rightarrow \Hom(\Coh_X, \Coh_S)$$ is fully faithful,
and its essential image consists of precisely those functors $F$
which extend to continuous, tame, additive tensor functors
$\hat{F}: \QC_X \rightarrow \Mod_{\calO_S}$.
\end{corollary}

\begin{remark}
It is unfortunate that there does not seem to be any simple
criterion on the functor $F$ which may be used to test whether or
not $\hat{F}$ is tame.
\end{remark}

We are now ready to give the proof of Theorem \ref{decoy}.

\begin{proof}
We have already noted that $\Hom(S^{\an}, X^{\an}) =
\Hom_{\C}(S^{\an}, X)$. Consequently, it will suffice to prove
that the natural functors $\Hom(S,X) \rightarrow \Hom(S^{\an},X)$
and $\Hom(S, \Spec \C) \rightarrow \Hom(S^{\an}, \Spec \C)$ are
equivalences. We will focus on the former (the latter is just the
special case where $X = \Spec \C$).

We are now free to drop the assumption that $X$ is of finite type
over $\C$ (or even that $X$ is an algebraic stack over $\C$); all
we will need to know is that $X$ is Noetherian and geometric. We
can therefore apply Corollary \ref{corevent} to deduce that
$\Hom(S^{\an}, X)$ and $\Hom(S,X)$ are equivalent to full
subcategories $\calC_0 \subseteq \calC = \Hom(\Coh_X,
\Coh_{S^{\an}})$ and $\calC'_0 \subseteq \calC' = \Hom(\Coh_X,
\Coh_S)$.

Using Serre's GAGA principle, we
may deduce that $\Coh_{S}$ is equivalent to $\Coh_{S^{\an}}$ as an
abelian tensor category. It follows that we may identify $\calC$
with $\calC'$. To complete the proof, it suffices to show that the
subcategories $\calC_0 \subseteq \calC$ and $\calC'_0 \subseteq
\calC'$ coincide (under the identification of $\calC$ with
$\calC'$). Let $F: \Coh_{X} \rightarrow \Coh_S$ be any additive
tensor functor. Then $F$ admits continuous extensions
$$ \hat{F}: \QC_X \rightarrow \Mod_{\calO_S}$$
$$ \hat{F}': \QC_X \rightarrow \Mod_{\calO_{S^{\an}}}.$$
We must show that $\hat{F}$ is tame if and only if $\hat{F}'$ is
tame.

The essential point is to observe that $\hat{F}' M = p^{\ast}
\hat{F} M$, where $p: (S^{\an}, \calO_{S^{\an}}) \rightarrow
(S,\calO_S)$ is the natural map. It follows immediately that if
$\hat{F}$ is tame, then $\hat{F}'$ is tame. For the converse, let
us suppose that $\hat{F}'$ is tame. Let $M \in \QC_{X}$ be flat.
Then $\hat{F}' M \in \Mod_{\calO_{S^{\an}}}$ is flat, so that the
stalk $(\hat{F}' M)_s = (\hat{F} M)_s \otimes_{\calO_{S,s}}
\calO_{S^{\an},s}$ is flat over $\calO_{S^{\an},s}$ at every
closed point $s$ of $S$. Since $\calO_{S^{\an},s}$ is faithfully
flat over $\calO_{S,s}$, we deduce that $(\hat{F} M)_s$ is flat
over $\calO_{S,s}$ at every closed point $s$ of $S$. Since
$\hat{F} M$ is a quasi-coherent sheaf on $S$, this implies that
$\hat{F} M$ is flat.

Let us now suppose that $$0 \rightarrow M' \rightarrow M
\rightarrow M'' \rightarrow 0$$ is a short exact sequence of
quasi-coherent sheaves on $X$, with $M''$ flat. We wish to show
that this sequence remains exact after applying the functor
$\hat{F}$. It suffices to show that the map $f: \hat{F} M'
\rightarrow \hat{F} M$ is injective. The kernel of $f$ is
quasi-coherent. Consequently, if the kernel of $f$ is nonzero,
then it has a nonzero stalk at some closed point $s \in S$. Since
$\calO_{S^{\an},s}$ is faithfully flat over $\calO_{S,s}$, we
deduce that the map $\hat{F}' M' \rightarrow \hat{F}' M$ is not
injective at the point $s \in S^{\an}$, a contradiction.
\end{proof}

\begin{remark}
The preceding argument is in fact quite general. It requires only the following facts:

\begin{enumerate}
\item The analytification $S^{\an}$ is local for the {\et}ale topology.
\item The map $(S^{\an},\calO_{S^{\an}}) \rightarrow (S,\calO_S)$ is faithfully flat.
\item Analytification induces an equivalence of categories from the category of coherent sheaves on $S$ to the category of coherent sheaves on $\calO_{S^{/an}}$.
\end{enumerate}

These hypotheses are frequently satisfied in other circumstances; for example, whenever $S$ is a Deligne-Mumford stack which is proper over $\C$. They are also satisfied if one replaces complex analytification by formal completion or rigid analytification (provided that one employs the correct {\et}ale topology on the ``analytic'' side).
\end{remark}

\begin{remark}
Theorem \ref{decoy} is not necessarily true if the stack $X$ is
not assumed to be geometric. For example, it can fail if $X$ is
the classifying stack of an abelian variety.
\end{remark}

\end{document}